\documentclass[11pt,psfig]{article}
\usepackage[dvips]{graphicx}
\usepackage{epsfig,amssymb,latexsym,amsmath}
\usepackage{euscript}
\usepackage{cite}
\usepackage{euscript}

\usepackage{cancel}
\usepackage{tikz}

\newtheorem{theorem}{Theorem}%[section]
\newtheorem{definition}{Definition}%[section]
\newcommand{\bea}{\begin{eqnarray}}
\newcommand{\eea}{\end{eqnarray}}
\newcommand{\be}{\begin{equation}}
\newcommand{\ee}{\end{equation}}
\newcommand{\ba}{\small\begin{array}}
\newcommand{\ea}{\end{array}\normalsize}

\newcommand{\G}{\EuScript{G}}

\newcommand{\rl}{{\bf R}^1}

\newcommand{\rr}{{\bf R}}
\newcommand{\rk}{{\bf R}^k}

\newcommand{\rn}{{\bf R}^n}
\newcommand{\rbn}{{\bf R}^N}

\newcommand{\nn}{{\bf N}}

\newcommand{\zz}{{\bf Z}}

\newcommand{\newdelta}{\delta}

\newcommand{\rtwo}{{\bf R}^2}

\newcommand{\rmiplusl}{{\bf R}^{1}}

\newcommand{\rmpplust}{{\bf R}^{1}}
\newcommand{\rkplusmpplusl}{{\bf R}^{k {+}1}}
\newcommand{\rkplusmpplust}{{\bf R}^{k {+} 2}}
\newcommand{\rml}{{\bf R}^{1}}
\newcommand{\rmt}{{\bf R}^{1}}

\newcommand{\rmpplusl}{ {\bf R}^{1} }

\newcommand{\rmi}{{\bf R}^{1}}

\begin{document}
\parindent=0pt
\parskip=.2cm
\markboth{ddf}{fdfd}

\begin{center}
 \Large \bf  Input-to-state stabilization of the perturbed  systems in the generalized triangular form
 \end{center}

\begin{center}
Sergey Dashkovskiy and Svyatoslav S. Pavlichkov
\end{center}
\begin{center}
Center of Industrial Mathematics, Faculty Mathematics and Computer Science,
        University of Bremen, Bremen, Germany
\end{center}
\begin{center}
Faculty Mathematics and Computer Science, Taurida National University, Simferopol, Ukraine
\end{center}

%%%%%%%%%%%%%%%%%%%%%%%%%%%%%%%%%%%%%%%%%ABSTRACT

{\bf Abstract:} We consider nonlinear control systems of the so-called generalized triangular form (GTF) with time-varying and periodic dynamics which linearly depends on some external disturbances.
Our purpose is to construct a feedback controller which provides the global input-to-state stability of the corresponding closed-loop w.r.t. the disturbances.
To do this, we combine the method proposed in the earlier work \cite{pavlichkov_ge_2009} devoted the the global asymptotic stabilization of the GTF systems without disturbances with the ISS theory
for time-varying systems proposed in \cite{wang}.
Following this pattern we construct a feedback which provides the properties of uniform global stability and asymptotic gain w.r.t the disturbances. Then we
obtain the semi-uniform ISS of the closed-loop system.


 \baselineskip=1.5\normalbaselineskip

%% main text

\section{INTRODUCTION}

 One of the most popular framework for design in nonlinear control theory
is backstepping. Originally, this approach was proposed for constructing
Lyapunov stabilizers \cite{kokotovic_sussmann,saberi,coron_praly}; very soon this technique was applied to
solving adaptive control problems: firstly when the dynamics of a strict-feedback system is linear
w.r.t. unknown parameters \cite{kanellakopoulos_kokotovic_morse_1991,Krstic_Kanellakopoulos_Kokotovic_1992_nachalo,krstic_kanellakopoulos_kokotovic_book,annaswami}, then these results were extended to the cases of nonlinear parametrization \cite{kojic_annaswamy_Automatica_2002}, unknown control directions \cite{jiang_1998}, backstepping for the systems with time delays \cite{ge_hong_lee_2003_TF_with_delays,ge_hong_lee_2005_TF_with_delays}, backstepping for the Volterra systems \cite{kps2} etc. 
 Let us remark that the classical version of this approach is applicable to the so-called strict-feedback form or, more generally, 
to the triangular form in the so-called regular case (the latter being intorduced in 1973 in \cite{kor1}), i.e., when the triangular system is feedback linearizable. As the exception we can mention works devoted to 
polynomial extensions of the strict-feedback forms \cite{celikovsky,tzamtzi_tsinias_99,wlin,qian_lin_2006_output_stabilization_poly_TF} as well as a more general situation \cite{tsinias_95,Cel_Nij}. This leads to the concept of the so-called generalized triangular form - \cite{pavl1,kp2,pavlichkov_ge_2009} (next called GTF). In the latter works the problem
of global robust controllability and that of global asymptotic stabilization of generalized triangular form systems was successively solved.

 On the  other hand, in many applications one has to consider systems subject to disturbances. In this case the input-to-state stability (ISS) framework introduced
in \cite{sontag} is very fruitful for stability analysis.
Therefore, having obtained the results on global asymptotic stabilization for the GTF systems \cite{pavlichkov_ge_2009}, it is natural to consider a GTF system with some external disturbances
in its dynamics and to ask whether it is possible to construct a feedback controller which provides a global input-to-state stability property with respect to the disturbances.
The goal of the current paper is to extend the result of the work \cite{pavlichkov_ge_2009} to this situation. 
In our case we will use the notion of uniform ISS developed for the case of time-varying systems \cite{wang}.
 A similar problem was considered for systems of the strict-feedback form in \cite{freeman_krstic_kokotovic_98} (however the strict-feedback form systems under consideration were not only with external disturbances but also with unknown parameters).  Since the GTF is an extension of the TF and strict-feedback forms, we extend the results
 of \cite{freeman_krstic_kokotovic_98} in the current paper in this sense as well.

\section{PRELIMINARIES}

Throughout the paper,  ${\bf N}$ and ${\bf Z}$ denote the
sets of all natural and integer numbers respectively,
$\langle \cdot,\cdot\rangle$  denote the scalar product in
$\rbn$ (for any $N \in \nn$); for $A \subset {\bf R}$, ${\rm mes} A$ and $\overline{A}$
denote the Lebesgue measure (if $A$ is
measurable) and the closure of $A$  respecitvely. For a vector $\xi{\in}\rbn,$ by
$|\xi|$ we denote its quadratic norm, i.e.,
$|\xi|{=}\langle\xi,\xi\rangle^{\frac{1}{2}}$.

A function $\alpha$ of ${\bf R}_{+}$ to ${\bf R}_{+}$ is said to be of class $\EuScript{N}$ if it is continuous and non-decreasing, is of class
$\EuScript{K}$ if it is continuous, positive definite and strictly increasing; and is of class $\EuScript{K}_{\infty}$ if it is of class $\EuScript{K}$ and unbounded.
A function $\beta$ of ${\bf R}_{+} \times {\bf R}_{+}$ to ${\bf R}_{+}$ is said to be of class  $\EuScript{KL}$ if for each fixed $t \ge 0$ the function $\beta(\cdot,t)$
is of class $\EuScript{K}_{\infty}$ and for each fixed $s \ge 0, $ we have $\beta(s,t) \rightarrow 0$ as $t \rightarrow +{\infty}$ and $t \mapsto \beta (s,t)$ is decreasing
Given any $\Delta(\cdot)$ in $L_{\infty}$ by $\parallel \Delta(\cdot) \parallel$  denote its $L_{\infty}$ - norm on $[0,+\infty[.$

Consider the nonlinear system
\begin{equation}
\dot x = f(t,x,\Delta) \label{sys_def}
\end{equation}
with states $x \in {\bf R}^n$ inputs $\Delta \in {\bf  R}^m, $ where $f$ is continuous w.r.t $(t,x,\Delta)$ and satisfies the local Lipschitz condition w.r.t. $(x,\Delta),$ whose
solution of the Cauchy problem $x(t_0)=x^0$ with $\Delta=\Delta(t)$ is denoted by $x(t,\xi,t_0,\Delta(\cdot)).$ Given any $\Delta(\cdot)$ in $L_{\infty}$ by $\parallel \Delta(\cdot) \parallel$  denote its $L_{\infty}$ - norm on $[0,+\infty[.$

The following three defintions and Theorem are borrowed from \cite{wang}
\begin{definition}
{ \em System (\ref{sys_def}) is input-to-state stable (ISS) iff there are $\beta \in \EuScript{KL},$ $\Upsilon_0 \in \EuScript{N}$ and $\gamma \in \EuScript{K}$ such that
for each $t_0,$ each $\xi$ and each $\Delta(\cdot),$ we obtain for all $t \ge t_0$
\begin{displaymath}
|x(t,\xi,t_0,\Delta(\cdot))| {\le} \beta(\Upsilon_0(t_0)|\xi|, t{-}t_0) {+} \gamma(\parallel \Delta(\cdot) {\parallel}_{L_{\infty} [t_0, +\infty[})
\end{displaymath}
System (\ref{sys_def}) is semi-uniformly input-to-state stable if it is ISS and furthermore there exists $\Upsilon (\cdot) \in \EuScript{K}$ such that
\begin{equation}
|x(t,\xi,t_0,\Delta(\cdot))| \le \max \{ \Upsilon (|\xi|), \Upsilon (\parallel \Delta(\cdot) \parallel) \} \; \; \forall t \ge t_0 \label{def_1}
\end{equation}
for all $\xi \in {\bf R}^n,$ and $\Delta(\cdot)$ in $L_{\infty}.$ }
\end{definition}

\begin{definition}
{\em 
We say that system (\ref{sys_def}) satisfies uniform local stability (ULS) property if there are  $\Upsilon \in \EuScript{K},$ and $\delta > 0$  such that
 for all $|\xi| \le \delta$ and all $\parallel \Delta(\cdot) \parallel \le \delta$ we have (\ref{def_1}). We say that system (\ref{sys_def}) satisfies uniform global stability (UGS) property
if $\delta=+\infty,$ i.e., (\ref{def_1}) holds for all $\xi$ and $\Delta(\cdot).$ }
\end{definition}

\begin{definition}
{\em
We say that system (\ref{sys_def}) satisfies the asymptotic gain (AG) property if there is $\gamma (\cdot )\in \EuScript{K}$ such that 
\begin{displaymath}
\limsup\limits_{t\rightarrow +\infty}|x(t,\xi,t_0,\Delta(\cdot))| \le \gamma(\parallel \Delta(\cdot) \parallel).
\end{displaymath}}
\end{definition}

\begin{theorem} \cite{wang}
 {\em System (\ref{sys_def}) is semi-uniform ISS if and only if it is ULS and AG.}
\end{theorem}

%The most of existing methods for stability investigation of nonlinear systems are based on Lyapunov functions or more generaly storage functions.
%Different method of stabilization were developed resently with Lyapunov approach: backstepping, forwarding, Zubov-method? ....
%%feedback linarization, flatness...
%
%Forwarding....
%
%Backsteping ...
%
%
%Framework of dissipative systems is suitable..... % some words about it
%
%Passive systems % this is a subcluss of dissip sys.
%
%In case of disturbances a very appropriate notion that was introduced by E. Sontag in \cite{Son89}.
%
%ISS systems % this is aniother subcluss of dissip sys.
%very goog + many results...

%%%%%%%%%%%%%%%%%%%%%%%%%%%%%%%%%%%%%%%%%%%%%%%%%%%%%%%%%%%%%%%%%%%%%%%%%%%%%%%%
\section{MAIN RESULT}

 %Also, we use the
%French notation for open intervals, i.e., for $a<b$ (with $a \in
%{\bf R}$ and $b \in {\bf R}$), $]a,b[$ denotes the open interval
%$]a,b[:=(a,b) = \{ x \in {\bf R} \;  | \; a<x<b \; \}.$

We consider the following system
\begin{equation}
 \left\{ \begin{array}{l}
    \dot x_1 = f_1(t,x_1,x_2) + {\newdelta}_1 (t) \Phi_1(t,x_1) \\
    \dot x_2= f_2(t,x_1,x_2,x_3)+ {\newdelta}_2 (t) \Phi_2(t,x_1,x_2)\\
     \ldots \\
\dot x_n= f_{n}(t,x_1,...,x_n,u) + {\newdelta}_n (t) \Phi_n(t,x_1,...,x_n) \end{array}\right.
 \label{eq.main.system.1}
\end{equation}
 where $u \in
\rl$ is the control, $x=[x_1,...,x_{n}{]}^T \in \rn$
is the state and ${\newdelta}_1 (t),{\newdelta}_2 (t),...,{\newdelta}_n (t)$ are some external disturbances
(in general ${\newdelta}_i (t)$ can be vectors of different finite dimensions).
We assume that (\ref{eq.main.system.1}) satisfies the following assumptions:

\begin{enumerate}

\item[A1:] $f=(f_1,...,f_n{)}^T$ and $\Phi_i$ are of class $C^{n{+}1}$
and $T$-periodic in time with some $T>0,$ i.e.,
$f(t+T,x,u)=f(t,x,u)$ and $\Phi_i(t+T,x)=\Phi_i(t,x)$  for all $[t,x,u]$ in $\rr\times\rn\times\rl.$

\item[A2:] $ f_i (t,x_1,...,x_i,\cdot): \rmiplusl \rightarrow
\rmi$ is a surjection, i.e., $ f_i (t,x_1,...,x_i,\rmiplusl){=}
\rl$ for
    every
 $[x_1,...,x_i] \in \rml \times \ldots\times \rmi,$ and every
 $t\in[0,T],$ $ i=
1,\ldots,n. $

\item[A3:]  there exist  $x_i^{\ast} \in  \rmi,$  $1\le i\le
n,$  and  $u^{\ast}=x_{n{+}1}^{\ast}$ in  $\rl$  such that
$
 \frac{\partial f_i}{\partial x_{i+1}} (t,x_1^{\ast},...,
x_{i+1}^{\ast}) \not= 0$  for every $t\in[0,T],$ $i=1,...,n,$  and
such that $f(t,x^{\ast},u^{\ast})=\Phi_i(t,x^{\ast})=0$ for all $t \in [0,T],$ $i=1,\ldots,n$

\end{enumerate}

The following example shows that even global asymptotic stabilization
of the time-invariant triangular systems is not always possible if one wants to use
a $C^1$ - feedback of the form $u=u(x)$. On the other hand, as we can see below, if we allow the feedback to be a time-varying,
it can reslove the problem (even a periodic feeback will suit). That is why we start with the
$T$ - periodic systems (of course, our result will be applicable to the time-invariant dynamics as a partial case as well).

  {\bf  Example 1.} \cite{kp2,pavlichkov_ge_2009} Consider the system
\begin{equation}\left\{ \begin{array}{l}
\dot x_1 = x_2^3 -  (1-x_1^2)x_2 \\ \dot x_2 = u
\end{array}\right. \label{eq.example.polynomial}
\end{equation}
and suppose there is a feedback  $u=u(x_1,x_2)$ of class
$C^1,$ which {\em globally } asymptotically stabilizes
(\ref{eq.example.polynomial}) into $[0,0].$ Put: $g(x):=
[x_2^3{-}(1{-}x_1^2)x_2, u(x_1,x_2)]^T,$ and  $C:= \{
[x_1,x_2] {\in} \rtwo \; | \; x_1^2+x_2^2 =1 \}.$ Since the feedback
$u=u(x)$ is continuous on $C,$ and  globally  stabilizes
(\ref{eq.example.polynomial}), we have $u(x)\not=0$ for all
$x\in C.$ Then, the map $C{\ni}x\mapsto \frac{g(x)}{|g(x)|} = [0,
\frac{u(x)}{|u(x)|}{]}^T$ is well-defined. On the one hand there is a homotopy between the map and  $C \ni x
\mapsto (-x) \in C $ (see the proof of the famous Brockett theorem
\cite{brockett_1983} given in \cite{sontag1}, p. 184), but on the other these maps have different degrees. This
contradiction proves that there is no a feedback
$u= u(x)$ of class $C^1$ which globally stabilizes
(\ref{eq.example.polynomial}).

%
%\begin{remark}\label{r2}
%On the other hand, it appears that for each $T>0,$ it is possible to
%construct a $T$ - periodic feedback of class $C^1$ which stabilizes
%the system (\ref{eq.example.polynomial}) globally. Furthermore, even
%if the triangular system is non-autonomous but $T$ - periodic with
%an arbitrary $T>0$ (and satisfies A1-A3 of course), then it is also
%possible to stabilize the system globally by means of a $T$ -
%periodic control of class $C^1$ which is the main result of our
%paper to be given below. Therefore, it is natural to formulate the
%stabilization problem for $T$ - periodic dynamics with some $T>0.$
%In Remark \ref{r4} (see Section \ref{section_four}), we discuss the
%difficulties that appear when trying to remove the assumption on
%T-periodicity.
%\end{remark}

%However, global stabilization of (\ref{eq.main.system.1}) is easy in
%some cases when (\ref{eq.main.system.1}) does not satisfy the
%above-mentioned Tsinias conditions A3 (i), (ii) from
%\cite{tsinias_95}. The following system can be presented as an
%example
% \begin{equation}\left\{ \begin{array}{l}
%\dot x_1 = g (x_2); \\ \dot x_2 = u  \end{array}\right.
%\label{eq.example.singular.backstepping}
%\end{equation}
%where $g(\cdot)$ of class $C^{3} (\rr;\rr)$ is given by
%\begin{equation*} g(x_2)= 0,  \; \;  \mbox{ if } \;  x_2
%< {-}2 \pi; \; \; \;  \mbox{ and }
%\end{equation*}
%\begin{equation*}
% g(x_2)=
% (x_2+2\pi)^5  {\sin}(x_2+2\pi), \; \;  \mbox{ if } \;  x_2 \ge
% {-}2\pi.
% \end{equation*}
%and its equilibrium point is $[0,0].$

Our main result is as follows.
\begin{theorem}\label{T1} {\em
Assume that system (\ref{eq.main.system.1}) satisfies
conditions A1-A3. Then, for any $\mu \in {\bf N} \cup \{ + \infty \}$ system (\ref{eq.main.system.1}) is globally
semi-uniformly input-to-state stabilizable into $x^{\ast}$ by means of a feedback law $u(t,x)$ of class
$C^{\mu}(\rr\times\rn;\rl)$ such that $u(t+T,x)=u(t,x)$ for all
$[t,x]\in\rr{\times}\rn$ and $u(t,x^{\ast})=u^{\ast}$ for all $t \in
\rr,$ where $T>0$ is the period mentioned above in A1.} %and such that % , if $\newdelta_i(t)=0,$ then 
%the closed-loop system with this feedback 
%is uniformly globally asymptotically stable (w.r.t. the equilibrium point $x^{\ast}$) and the asymptotic gain property holds (w.r.t. the equilibrium point $x^{\ast}$).
%and, for each $t_0 \in \rr$  if $\sum\limits_{i=1}^{n} \parallel \delta_i(\cdot) {\parallel}_{L_{\infty} [t_0, t_0+T]}^2 < d $ for all $T>0,$
%then any trajectory $x(t)$ of the system (\ref{eq.main.system.1}) with
%$u=u(t,x),$ defined by $x(t_0)=x_0$ satisfies $\limsup\limits_{t \rightarrow +\infty} |x(t)| \le \gamma(d),$ whatever $x_0 \in \rn$.
\end{theorem}

Let us remark that $x_i$ and $u$ can be vectors in general as in \cite{pavlichkov_ge_2009}
and we assume them to be scalar for the simplicity only (for vectors, the argument will be similar)

%%%%%%%%%%%%%%%%%%%%%%%%%%%%%%%%%%%%%%%%%%%%%%%%%%%%%%%%%%%%%%%%%%%%%%%%%%%%%%%%
\section{BACKSTEPPING DESIGN}

Let $k$ be in $\{0,...,n{-}1 \}.$ %Define $k:=m_1+ \ldots +m_p, $
%if $p\ge1,$ and $k=0$ if $p=0.$
For each $y_0\in \rkplusmpplusl,$ each
$\omega_0\in{\rmpplusl},$ and each $r>0,$ let $B_r(y_0)$ and ${\Omega}_r(\omega_0)$ denote the open balls
 \begin{equation*}
B_r(y_0):= \{ y\in {\rkplusmpplusl} \; | \; \; |y-y_0|<r  \}; \; \; \; 
 {\Omega}_r({\omega}_0):= \{ \omega\in \rmpplusl \; | \; \;
|\omega-\omega_0|<r \}
\end{equation*}
and ${\overline B}_r(y_0)$ and ${\overline \Omega}_r(\omega_0)$ be
their closures.

Consider a control system
\begin{equation}
\dot z = g(t,z,z_{k{+}1}) + \sum\limits_{j=1}^{N_k} {\Delta}_{j}(t) \varphi_j (t,z),  \; \; \; t\in\rr
\label{eq.backstepping.subsystem.1}
\end{equation}
where $z_{k{+}1} \in \rmpplusl$ is the control,
$z=[z_1,...,z_k{]}^T \in \rk,$ is the state, and $\delta (t) = [{\Delta}_1(t),...,{\Delta}_{N_k} (t)]$ is some external disturbance.

Following \cite{pavlichkov_ge_2009},  we also consider
a dynamical extension of (\ref{eq.backstepping.subsystem.1}),  i.e., the system
\begin{equation*}
\left \{ \begin{array}{l}
\dot z = g(t,z,z_{k{+}1}) + \sum\limits_{j=1}^{N_k} {\Delta}_{j}(t) \varphi_j (t,z) \\
\dot z_{k{+}1} {=} g_{k{+}1}(t,z,z_{k{+}1},v) {+} \sum\limits_{j=1}^{N_k} {\Delta}_{j}(t) {\varphi}_{k{+}1,j} (t,z)
{+}\sum\limits_{j=N_{k}+1}^{N_{k+1}} {\Delta}_{j}(t) {\varphi}_{k{+}1,j} (t,z,z_{k{+}1})
  \end{array} \right.
  \end{equation*}
which we rewrite in the following vector form
\begin{equation}
\dot y = \psi (t,y,v) + \Delta(t) \phi(t,y),  \; \; \; \; \; \; t\in\rr,
\label{eq.backstepping.extended.subsystem.2}
\end{equation}
where $y=[z,z_{k{+}1}{]}^T \in \rkplusmpplusl$ is the state,
$v \in \rmpplust$ is the control, $\Delta(t)=[{\Delta}_1(t),...,{\Delta}_{N_k} (t), {\Delta}_{N_{k}+1} (t),...,{\Delta}_{N_{k+1}} (t) ]$ is its external disturbance (with $N_{k+1}>N_k$), and $\psi(t,y,v)$  and
$\phi(t,y)$ are
given by

\begin{equation*}
\psi(t,y,v) = \left[ \begin{array}{l}g(t,y) \\
 g_{k{+}1}(t,y,v) \end{array} \right] \; \;
  \mbox{ and } \; \;  \phi (t,z,z_{k+1}) =
  \end{equation*}
  \begin{equation*}
  \left[\begin{array}{llllll} \varphi_1 & ... &  \varphi_{N_k} & 0 &...&0 \\
\varphi_{{k+1},1} & ...&  {\varphi}_{k+1, N_k} & {\varphi}_{k+1, N_k+1}  &...& {\varphi}_{k+1, N_{k+1}} \end{array} \right]
 \end{equation*}
 \begin{equation}
  \mbox{ for all }
   \; \; \;  [t,y,v]\in\rr\times \rkplusmpplusl.
\times\rmpplust \label{eq.backstepping.extended.subsystem.3}
\end{equation}
%with $g_{k{+}1} {\in} \rmpplusl.$

As in \cite{pavlichkov_ge_2009}, if $k{=}0,$ and system
(\ref{eq.backstepping.subsystem.1}) consists of $0$ equations, we
define $y:=z_{k{+}1}=z_1;$ $\psi(t,y,v):= g_{k{+}1}(t,y,v)=
g_1(t,z_1,v)$ with $v\in\rmt$ and we say that
(\ref{eq.backstepping.subsystem.1})  is empty or trivial and that
$\dot z_1 = g_1 (t,z_1,z_2)$ with states $z_1=y$ and controls
$z_2=v$ is the extension of the empty system
(\ref{eq.backstepping.subsystem.1}).

We assume that $\psi$ and $\phi$ satisfy the following
Assumptions:

\begin{enumerate}

\item[A1':] {\em  Functions $\psi$ and $\phi$ are of classes
$C^2(\rr\times\rkplusmpplust;\rkplusmpplusl)$ and $C^2(\rr\times\rkplusmpplusl;\rkplusmpplusl)$ respectively and there exists
$T>0$ such that $\psi(t+T,y,v)=\psi(t,y,v)$  and $\phi(t+T,y)=\phi(t,y)$ for all $[t,y,v]$ in
$\rr\times\rkplusmpplust$. }

\item[A2':] {\em For every $t\in\rr,$ we have:  $\psi(t,0,0)= 0 ;
$ and $  \frac{\partial g_{k{+}1}}{\partial v} (t,0,0) \not= 0
$.}

\item[A3':] {\em $g_{k{+}1} (t,y,\rmpplust) =\rmpplusl $ for every
$[t,y]\in[0,T]{\times}\rkplusmpplusl.$  }

\end{enumerate}

Given an initial state $z_0\in\rk,$ a feedback control
$\omega(t,z)$ of $\rr\times\rk$ to $\rmpplusl$ a disturbance $\delta(\cdot)$ and  $t_0\in\rr,$ let
$t \mapsto z(t,t_0,z_0,\omega(\cdot,\cdot), \delta(\cdot))$ denote the trajectory, of
system (\ref{eq.backstepping.subsystem.1}) that is defined by this
control $\omega(\cdot,\cdot),$ by this disturbance $\delta (\cdot)$ and by the initial condition
$z(t_0)=z_0.$ Similarly, for system
(\ref{eq.backstepping.extended.subsystem.2}), given an initial
state $y_0\in\rkplusmpplusl,$ a feedback $v(t,y)$ of
$\rr{\times}\rkplusmpplusl$ to $\rmpplust,$ a disturbance $\Delta(\cdot)$ and $t_0{\in}\rr,$ let
$y(t,t_0,y_0,v(\cdot,\cdot), \Delta(\cdot))$ denote the trajectory, of
(\ref{eq.backstepping.extended.subsystem.2}), that is defined by
the control $v(\cdot,\cdot),$ by the disturbance $\Delta(\cdot),$  and by the initial condition
$y(t_0)=y_0.$ %(If the controls $\omega$ and $v$ are not
%closed-loop but open-loop, or even constant, we use the same
%notation $z(\cdot,t_0,z_0,\omega),$  and $y(\cdot,t_0,y_0,v)$ as
%well.
 In addition, we presume that the existence and the
uniqueness of the solution of the corresponding Cauchy problem are
ensured in this definition. 
Of course, if $\omega$ and $v$ are at least of class $C^1,$  and if the disturbances are of class $L_{\infty},$ then
it guarantees the existence and the uniqueness of the
corresponding solution automatically.

Following \cite{pavlichkov_ge_2009},  for systems (\ref{eq.backstepping.subsystem.1}) and
(\ref{eq.backstepping.extended.subsystem.2}), we consider the
following Lyapunov pairs:
\begin{equation*}
V_k(z):=\langle z,z \rangle,  \; V_{k{+}1}(y):=\langle y,y
\rangle 
  \mbox{ for all }  z{\in}\rk;\; y{\in}\rkplusmpplusl
\end{equation*}

We reduce Theorem 2 to the following Theorem.

 \begin{theorem}\label{T2} {\em
Assume that systems
(\ref{eq.backstepping.subsystem.1}) and
(\ref{eq.backstepping.extended.subsystem.2}) satisfy
Assumptions A1'-A3'. Suppose there exist  sequences $\{ r_q
{\}}_{q{=}2}^{{+}\infty} {\subset} \rr,$  $\{ {\rho}_q
{\}}_{q{=}1}^{{+}\infty} {\subset} \rr$ and $\{ d_q
{\}}_{q{=}1}^{{+}\infty} {\subset} \rr$ such that $0<{\rho}_q<
r_{q{+}1}< {\rho}_{q{+}1}$ and  $0<d_q<d_{q+1}$ for all $q{\in}\nn$ such that
$r_q{\rightarrow}{+}\infty,$ ${\rho}_q{\rightarrow}{+}\infty$ and $d_q{\rightarrow}{+}\infty$ as
$q{\rightarrow}\infty.$ Assume that there exists a function ${\gamma}_{k} (\cdot) \in {\EuScript{K}}_{\infty}$ such that
$d_1 < \max\limits_{|z|\le {\rho}_1}  V_k(z)$ and the following conditions hold:}
\begin{enumerate}
\item[C1:] $\frac{\partial V_k (z)}{\partial z} (g(t,z,0) +  \sum\limits_{j=1}^{N_k} {\Delta}_{j} \varphi_j (t,z))  \le -
V_k(z) + {\gamma}_{k} (|\delta|)$ for all $\delta \in {\bf R}^{N_k},$ whenever $|z|^2<r_2^2,$ $z\in\rk,$ $t\in [0,T]$.
\item[C2:] For every $z_0{\in}\rk,$ and every  $t_0{\in} [0,T]$ if
$|z_0|^2{\le}r_{q{+}2}^2$ with some $q\in\nn$ then
\begin{equation*}
|z(t,t_0,z^0,0,\delta(\cdot))|^2\le \rho_{q{+}2}^2 - \frac{t-t_0}{T}
(\rho_{q{+}2}^2-\rho_q^2) \mbox{ for all }   t \in [t_0,t_0{+}T],
\end{equation*}
 whenever $\delta (\cdot)$ in $L_{\infty} [t_0,t_0+T]$ satisfies
${\gamma_k(\parallel \delta(\cdot) {\parallel}_{L_{\infty} [t_0,t_0+T]})} {\le} d_q.$
\end{enumerate} {\em
Then, for every $\mu\in {\nn}\cup \{\infty\},$  there exist
$q_0{\ge}0$ ($q_0{\in}{\bf Z}$) positive real numbers $r_1,$
$r_0,$ ..., $r_{-q_0},$  postitve real $d_0,$
$d_{-1},$ ..., $d_{{-}q_0{-}1},$ a sequence of positive real numbers  $\{
R_q {\}}_{q{=}{-}q_0{-}1}^{\infty}, $ a function ${\gamma}_{k+1} (\cdot) \in {\EuScript{K}}_{\infty}$ such that ${\gamma}_{k+1} (|\Delta|) \ge  {\gamma}_{k} (|\delta|)+ |\Delta|^2$ and a feedback controller
$v(\cdot,\cdot)$ of class
$C^{\mu}(\rr\times\rkplusmpplusl;\rmpplust)$   such that $0<{R}_q<
r_{q{+}1}< {R}_{q{+}1}$ and $0<d_q<d_{q+1}$ for all $q\ge{-}q_0-1,$ $q\in{\bf Z}$ and $d_{{-}q_0{-}1} < \max\limits_{|y| \le R_{{-}q_0{-}1}}  V_{k+1}(y)$ and
such that the following conditions hold:} 
\begin{enumerate}
\item[(i)] $v(T{+}t,y){=}v(t,y) $ for all $[t,y]$ in
$\rr{\times\rkplusmpplusl},$ and $v(t,0){=}0{\in}\rmpplust$ for
all $t{\in}\rr.$
\item[(ii)] For each $t{\in}\rr,$  each $y{=}
[z,z_{k{+}1}{]}^T {\in}  {\overline B}_{r_{{-}q_0}}(0),$ and each $\Delta {\in} {\bf R}^{N_{k{+}1}},$ we have:
\begin{equation*}
\frac{\partial V_{k{+}1} (y)}{\partial y} (\psi(t,y,v(t,y)) {+} \Delta \phi(t,y))  {\le} {-}
V_{k{+}1}(y)
{+} {\gamma}_{k{+}1} (|\Delta|)   
\end{equation*}
\item[(iii)] For every $y_0{\in}\rkplusmpplusl,$ and every
$t_0\in {\bf R}$ if $|y_0|^2\le r_{q{+}2}^2$ with some $q\ge - q_0
- 1,$ $q \in {\bf Z}$ then
\begin{equation*}
|y(t,t_0,y^0,v(\cdot,\cdot), \Delta(\cdot))|^2\le R_{q{+}2}^2 {-} \frac{t{-}t_0}{T}
(R_{q{+}2}^2{-}R_q^2) \mbox{ for all } t\in[t_0,t_0{+}T],
\end{equation*}
whenever $\Delta(\cdot) \in L_{\infty} [t_0,t_0+T]$ satisfies $\gamma_{k{+}1}(\parallel \Delta(\cdot) {\parallel}_{L_{\infty} [t_0,t_0+T]}) \le d_q.$
%\end{equation}
\end{enumerate}
{\em
(If $k{=}0,$ i.e., system (\ref{eq.backstepping.subsystem.1}) is
empty, we say that Conditions C1, C2 hold by definition, and the
Theorem states that, for the corresponding extension
(\ref{eq.backstepping.extended.subsystem.2}), there is a control
$v(\cdot,\cdot)$ such that Conditions (i), (ii), (iii) hold with $\gamma_1(|\Delta|) = |\Delta|^2$).}
\end{theorem}

It is easy to prove that Theorem 3 implies Theorem 2. Indeed, assume that system (\ref{eq.main.system.1}) satisfies Assumptions A1-A3 and 
without loss of generality assume that $x^{\ast}=0,$ $u^{\ast}=0.$ 

For $k{=}0,$ define $y:=x_1,$ $v:=x_2,$ $\psi:=f_1(t,y,v),$ $\phi := \Phi_1 (t,y),$ $\Delta := \delta_1$
and find the feedback $\alpha_1(t,y):= v(t,y),$ the $\EuScript{K}_{\infty}$  - function $\gamma_1(|\Delta|):= |\Delta|^2$
and positive numbers $r_q$ $(q \ge - q_0)$ and $R_q,$ 
$d_q$ $(q \ge - q_0 -1)$ satisfying all the statement of Theorem 3 including (i)-(iii). Without loss of generality, assume that $q_0=2$ (otherwise
shift the indexation accordingly). 

Then, for $k=1$ redefine:
\begin{displaymath}
z=z_1:=x_1,  \; \; \; z_{k+1}=z_2:=x_2-\alpha_1(t,x_1), 
\; \; \; \;  y:=[z_1,z_2], \; \; \; v:=x_3
\end{displaymath}
\begin{displaymath}
g(t,z,z_{k+1}) := f_1(t,z_1,z_2+\alpha_1(t,z_1));
\end{displaymath}
\begin{displaymath}
[\Delta_{1},...,\Delta_{N_{k}}]:=\delta_1, \; \; \; [\Delta_{N_{k}+1},...,\Delta_{N_{k+1}} ]:=\delta_2, 
\; \; \; [\varphi_{1},...,\varphi_{N_{k}}] (t,z) := \Phi_1(t,z), \; \; \; 
\end{displaymath}
\begin{displaymath}
g_{k+1}(t,z,z_{k+1},v):= f_2 (t,z,z_2+\alpha_1(t,z_1),v)
 - \frac{\partial \alpha_1 (t,z)}{\partial t} - \frac{\partial \alpha_1 (t,z)}{\partial z} g(t,z,z_2),
\end{displaymath}
\begin{displaymath}
[\varphi_{k+1,1},...,\varphi_{k+1,N_{k}}] (t,z):= - \frac{\partial \alpha_1 (t,z)}{\partial z} \Phi_1(t,z) 
\end{displaymath}
\begin{displaymath}
[\varphi_{k+1,N_{k}+1},...,\varphi_{k+1,N_{k+1}}] (t,z,z_{k+1}) :=
 {\Phi}_2(t,z_1,z_2+\alpha_1(t,z_1))
\end{displaymath}
Then, for $k{=}1,$ system (\ref{eq.backstepping.subsystem.1}) satisfies Assumptions C1-C2 of Theorem 3 and, applying Theorem 3, we obtain the existence of
$r_q,$ $R_q,$ $d_q,$ $v(t,y),$ and $\gamma_2(\cdot) \in {\EuScript{K}}_{\infty}$ satisfying the statement of Theorem 3 for $k=1$ (and satisfying (i)-(iii)).
Similarly, after $n$ coordinate transformations and $n$ steps of the backstepping procedure, we obtain  system (\ref{eq.backstepping.extended.subsystem.2}) of dimension $k+1=n$
and the existence of the corresponding $T$ -periodic feedback $v(t,y)$ of ${\bf R} \times {\bf R}^n$ to ${\bf R},$ $\gamma_n(\cdot) \in {\EuScript{K}}_{\infty}$ and the existence of positive numbers $r_q$ $(q \ge - q_0)$ and $R_q,$ 
$d_q$ $(q \ge - q_0 -1)$ satisfying (i)-(iii) and the statement of Theorem 3. 
Then one proves that the $n$ - dimensional closed-loop system
\begin{displaymath}
\dot y = \psi (t,y,v(t,y)) + \Delta(t) \phi(t,y),  \; \; \; \; \; \; t\in\rr,
\end{displaymath}
satisfies the asymptotic gain (AG) property
and the global unfiorm stability (UGS) property.

(AG): Take any  $\{ d_q {\}}_{q={-}q_0{-}2}^{{-}\infty}$  such that $0<d_q<d_{q{+}1}$ for all $q \in \zz$ and such that $d_q \rightarrow 0$ as $q \rightarrow \ {-}\infty.$ From (iii) we obtain:
\begin{displaymath}
\gamma_n (\parallel \Delta(\cdot) \parallel) \le d_q \Rightarrow
\limsup\limits_{t \rightarrow +\infty}|y(t,t_0,y^0,v(\cdot,\cdot), \Delta(\cdot))|^2 \le R_q^2, 
\end{displaymath}
\begin{equation}
\mbox{ whenever } \;  q\ge -q_0-1
\label{new_eq_1}
\end{equation}
and from (ii), (iii) we obtain:
\begin{displaymath}
\gamma_n (\parallel \Delta(\cdot) \parallel) \le d_q \Rightarrow 
\limsup\limits_{t \rightarrow {+}\infty}|y(t,t_0,y^0,v(\cdot,\cdot), \Delta(\cdot))|^2 \le d_q {<}d_{{-}q_0{-}1}{<}R_{{-}q_0{-}1}^2,
\end{displaymath}
\begin{equation}
 \; \mbox{ whenever } q < -q_0-1
\label{new_eq_2}
\end{equation}
for all $y^0 \in \rn.$ Find any $\gamma(\cdot) \in {\EuScript{K}}_{\infty}$ such that
\begin{equation}
d_q<\gamma_n(|\Delta|)\le d_{q+1}\;  \Rightarrow \;  \gamma (|\Delta|) \ge R_{q{+}1}^2 
\; \; \mbox{ for each } \;  q \ge -q_0-2
\label{new_eq_3}
\end{equation}
\begin{equation}
 d_q<\gamma_n(|\Delta|)\le d_{q+1} \;  \Rightarrow \;  \gamma (|\Delta|) \ge d_{q{+}2}
  \; \;  \mbox{ for each } \;  q < -q_0-2
\label{new_eq_4}
\end{equation}
Then from (\ref{new_eq_1})-(\ref{new_eq_4}) we obtain the AG property:
\begin{displaymath}
\limsup\limits_{t \rightarrow +\infty}|y(t,t_0,y^0,v(\cdot,\cdot), \Delta(\cdot))|^2 \le \gamma (\parallel \Delta (\cdot) \parallel)
\end{displaymath}
whatever $y^0 \in \rn.$

(UGS):  Take any  $\{ d_q {\}}_{q={-}q_0{-}2}^{{-}\infty}$  such that $0<d_q<d_{q{+}1}$ for all $q \in \zz$ and such that $d_q \rightarrow 0$ as $q \rightarrow \ {-}\infty.$
Also take $\{ R_q {\}}_{q={-}q_0{-}2}^{{-}\infty}$  such that $0<R_q<R_{q{+}1}$ for all $q \in \zz$ and  $R_q \rightarrow 0$ as $q \rightarrow \ {-}\infty$ and
such that $d_q{<} R_q^2$ for all $q \in \zz$ (note that for $q \ge -q_0-1$ the latter follows from Theorem 3).
>From (iii), we obtain that for each  $q \ge -q_0-1,$ if  $|y_0|^2 \le r_{q{+}2}^2$  and  $\gamma_n (\parallel \Delta (\cdot) \parallel) \le d_q$ then 
\begin{displaymath}
|y(t,t_0,y^0,v(\cdot,\cdot), \Delta(\cdot))|^2 \le R_{q+2}^2 \; \; \forall t \ge t_0.
\end{displaymath}
Similarly, from (ii) and from the inequalities $d_q{<} R_q^2{<}R_{q+2}^2,$ which hold true for all $q \in \zz,$
we obtain that for each  $q < -q_0-1,$ if  $|y_0|^2 \le r_{q{+}2}^2$  and  $\gamma_n (\parallel \Delta (\cdot) \parallel) \le d_q$ then 
\begin{displaymath}
|y(t,t_0,y^0,v(\cdot,\cdot), \Delta(\cdot))|^2 \le R_{q+2}^2 \; \; \forall t \ge t_0.
\end{displaymath}
Combining these two implications, we obtain the following one:
\begin{displaymath}
\forall \{ \overline{q}, \hat{q} \} \subset \zz \; \; (|y_0|^2 \le r_{{\overline{q}}{+}2}^2) \wedge (\gamma_n (\parallel \Delta (\cdot) \parallel) \le d_{\hat{q}}) 
\end{displaymath}
\begin{equation}
\Rightarrow \; |y(t,t_0,y^0,v(\cdot,\cdot), \Delta(\cdot))|^2 \le \max \{ R_{\overline{q}}^2, R_{\hat{q}}^2, \}
\label{new_eq_5}
\end{equation}
Find any $\Upsilon (\cdot) \in {\EuScript{K}}_{\infty}$ such that
\begin{displaymath}
r_{q+1}<|y_0|\le r_{q{+}2} \; \Rightarrow \; \Upsilon (|y_0|) \ge R_{q+2}^2
\end{displaymath}
\begin{displaymath}
d_{q-1}< \gamma_n (\parallel \Delta (\cdot) \parallel) \le d_{q} \; \Rightarrow \; \Upsilon (\parallel \Delta (\cdot) \parallel) \ge R_{q+2}^2
\end{displaymath}
for all $q \in \zz.$
Combining the latter with (\ref{new_eq_5}), we obtain the UGS:
\begin{displaymath}
|y(t,t_0,y^0,v(\cdot,\cdot), \Delta(\cdot))|^2 \le \max \{ \Upsilon (|y_0|), \Upsilon (\parallel \Delta (\cdot) \parallel)\}  
\end{displaymath}
for all $t \ge t_0,$ $y^0 \in \rn$

Since our transformation of coordinates was triangular, $T$ - periodic, and is a global diffeomorphsm of states, we see that the original system
(\ref{eq.main.system.1}) will also be UGS and AG with this feedback.
The proof of Theorem 2 is complete, it remains to prove Theorem 3.

%%%%%%%%%%%%%%%%%%%%%%%%%%%%%%%%%%%%%%%%%%%%%%%%%%%%%%%%%%%%%%%%%%%%%%%%%%%%%%%%

\section{PROOF OF THEOREM 3}

Following \cite{pavlichkov_ge_2009}, we prove the existence of numbers $r {\in}  ]0,\rho_1[,$ $d_{{-}q_0{-}1}$ in $]0,d_1[,$ feedback
$\nu(\cdot,\cdot)$ of class $C^{\infty}(\rr{\times}{\overline
B}_{2r}(0);\rmpplust)$ and function ${\gamma}_{k+1} (\cdot)$ of class ${\EuScript{K}}_{\infty}$ such that ${\gamma}_{k+1} (|\Delta|) \ge  {\gamma}_{k} (|\delta|)+ |\Delta|^2$ and $d_{{-}q_0{-}1} \le \max\limits_{|y|\le r}  V_{k+1}(y)$ and such that
\begin{equation}
\nu(t,0)=0;\;  \nu(t+T,y)=\nu(t,y)  \; \mbox{ for all } t
\in \rr,\;  y \in \rkplusmpplusl \label{eq.proof.beta.periodic.1}
\end{equation}
and 
\begin{equation*}
\frac{\partial V_{k{+}1} (y)}{\partial y} (\psi (t,y,\nu(t,y)) + \Delta \phi(t,y))
\le - V_{k{+}1} (y)
\end{equation*}
\begin{equation}
+ \gamma_{k+1} (|\Delta|)  \mbox{ for all }  \Delta {\in} {\bf R}^{N_{k{+}1}}, \;
y {=} [z,z_{k{+}1}] {\in}  {\overline B}_{2r}(0),\;  t{\in }\rr
\label{eq.proof.Lyapunov.inequality.2}
\end{equation}
Indeed, by condition C1 of Theorem 3, the derivative of
$V_{k{+}1}$ along the trajectories of
(\ref{eq.backstepping.extended.subsystem.2}) is
\bea \frac{dV_{k{+}1}}{dt} &=& \frac{\partial V_{k{+}1}}{\partial
y} (\psi (t,y,v) + \Delta \phi(t,y)) \nonumber\\ &=&  \frac{\partial V_{k}(z)}{\partial z}
(g(t,z,0) {+} \sum\limits_{j=1}^{N_k} {\Delta}_{j} {\varphi}_j (t,z))
{+}  \frac{\partial V_{k}(z)}{\partial z} \nonumber\\
&& (g(t,z,z_{k{+}1})  {-} g(t,z,0))  {+} 2 z_{k{+}1} (
g_{k{+}1}(t,y,v)  \nonumber\\
&+& \sum\limits_{j=1}^{N_{k{+}1}} {\Delta}_{j} {\varphi}_{k{+}1,j}(t,z,z_{k+1}))
{\le} {-} V_k(z) {+} \gamma_k (|\delta|) \nonumber\\
&+& z_{k+1} (2  g_{k+1} (t,y,v) {+} \frac{\partial V_k(z)}{\partial z}
J (t,z,z_{k+1}) \nonumber\\
&+&  2 \sum\limits_{j=1}^{N_{k{+}1}} {\Delta}_j \varphi_{k+1,j} (t,z,z_{k+1})) {\le} {-}V_k (z) + \gamma_k (|\delta|) + |\Delta|^2 \nonumber\\
 &+&  z_{k+1} (2 g_{k+1} (t,y,v)   
+  \frac{\partial V_k(z)}{\partial z}
J (t,z,z_{k+1})  \nonumber\\
&+&  z_{k+1} \sum\limits_{j=1}^{N_{k{+}1}} {\varphi}_{k{+}1,j}^2(t,z,z_{k+1})  )\nonumber
\eea
for all $\Delta \in {\bf R}^{N_{k+1}},$ whenever $ |z|^2 < r_2^2,$ $z \in \rk,$ $t \in [0,T],$
where
\begin{equation*}
 J (t,z,z_{k+1}) = \int\limits_{0}^{1} \frac{\partial g(t,z,\theta
z_{k{+}1})}{\partial z_{k{+}1}}d\theta
\end{equation*}

Then %%%%%%%using Lemma A1 from the appendix,
we obtain the existence of $r{\in}]0,\rho_1[$ and {\em $T-$
periodic} feedback $\nu(\cdot,\cdot)$ in $C^{\infty} (\rr{\times}{\overline
B}_{2r}(0); \rmpplust)$ such that
\begin{equation*}
 z_{k+1} (2 g_{k+1} (t,y,\nu(t,y))
+ \frac{\partial V_k(z)}{\partial z}
J (t,z,z_{k+1}) + z_{k+1}\sum\limits_{j=1}^{N_{k{+}1}} {\varphi}_{k{+}1,j}^2(t,z,z_{k+1}) )
\end{equation*}
\begin{equation*}
 \le
 {-} |z_{k{+}1}|^2 \; \; \; \;  \mbox{ for all }
  [t,y] \in  \rr \times  {\overline B}_{2r}(0).
 \end{equation*}
Take any $d_{{-}q_0{-}1}$ in $]0, d_1[ $  that satisfiy the inequality  $d_{-q_0-1} {<} \frac{1}{6} \max\limits_{|y| \le r} V_{k+1} (y)$ and any function ${\gamma}_{k+1} (\cdot)$ of class ${\EuScript{K}}_{\infty}$ such that ${\gamma}_{k+1} (|\Delta|) \ge  {\gamma}_{k} (|\delta|)+ |\Delta|^2.$ Then
(\ref{eq.proof.beta.periodic.1}),
(\ref{eq.proof.Lyapunov.inequality.2}) are satisfied, and
$\nu(\cdot,\cdot)$
 satisfies Condition (ii) of Theorem 3.

Let us point out that for $k=0$ all these arguments will be simplified (the terms corresponding to $z,$ $g(t,z,z_{k+1})$ and to their scalar product will be abscent) - similar remark can be made
for the next steps.

 Next, we extent our control onto the whole
state space to satisfy condition (iii).

 Define 
 \begin{displaymath}
 \varkappa :=\min \left\{ \frac{1}{6}\min\limits_{|z|{\ge}r}V_k(z), \;  \frac{1}{4}\left(\max\limits_{|z|{\le} \rho_1 }V_k(z) {-} d_1 \right)\right\}
 \end{displaymath}

Using the Gronwall-Bellman lemma and Condition C2 of Theorem 3,
%and (\ref{eq.proof.definition.of.sigma.3}),
%(\ref{eq.proof.definition.of.r_i.5}),
we find positive numbers $R_q > 0,$ $- q_0- 1
\le q \le 3,$ $\sigma_{-q_0}=\sigma_{-q_0+1}=\ldots
=\sigma_{1}= \sigma_{2}= \sigma_{3}=\sigma$ and $d_0>
d_{-1}> \ldots > d_{{-}q_0}$ with  $d_{{-}q_0} > d_{{-}q_0{-}1}$ (where $d_{{-}q_0{-}1}$ was chosen above) such that
first,
\begin{equation*}
\frac{\partial V_k (z)}{\partial z}( g(t,z,z_{k{+}1}) +  \sum\limits_{j=1}^{N_k} {\Delta}_{j} \varphi_j (t,z))
\end{equation*}
\begin{equation*}
= 2 \langle z,
g(t,z,z_{k+1}) +  \sum\limits_{j=1}^{N_k} {\Delta}_{j} \varphi_j (t,z) \rangle \le - 2 \varkappa, \mbox{ whenever }
\end{equation*}
\begin{equation*}
|z|^2{+} |z_{k{+}1}|^2<R_q \mbox{ and } \gamma_k(|\delta|) < d_q \; \; \; \;  \mbox{ for
all } \end{equation*}
\begin{equation*}
   y=[z,z_{k{+}1}]
  \in   \left( {\overline B}_{r_2}(0){\setminus} B_{R_{{-}q_0{-}1}}(0)\right)
  {\cap}
\left( \rk{\times}{\overline \Omega}_{3\sigma}(0) \right),
\end{equation*}
\begin{equation}
 - q_0- 1
\le q \le 1, \; \; \; t\in \rr
\label{eq.proof.definition.of.sigma.5_a}
\end{equation}
\begin{equation}
 r_q < \rho_q < R_q < r_{q{+}1}\; \; \; \; \; \; \;
\; \mbox{ for all } \; \; q = 1,2,3;
\label{eq.proof.conditions.for.R_i.6}
\end{equation}
\begin{equation*}  R_{-q_0}< 2r = 2r_{-q_0}, \; r_{-q_0}=r,  \mbox{ and  }
R_{q-1} < r_q < R_q < r_{q{+}1}
\end{equation*}
\begin{equation}
     \mbox{ for all }
-q_0\le q\le 1 \label{eq.proof.conditions.for.R_i.7}
\end{equation}
 second, for every $z_0\in \rk,$  every $t_0 \in [0,T],$
every $\omega(\cdot)$ in $C([t_0,t_0{+}T];\rmpplusl),$ and every $\delta(\cdot)$ in $L_{\infty} [t_0,t_0+T]$
\begin{equation*}
\mbox{ if } \;   \max\limits_{t_0 \le t \le
t_0{+}T}|\omega(t)| \le 3 \sigma_3,  \; \; |z_0|^2 \le r_{3}^2,
\;   \mbox{ and }  \; \gamma_k(\parallel \delta(\cdot) {\parallel}_{L_{\infty} [t_0,t_0{+}T]}) \le d_1,
\mbox{ then } 
\end{equation*}
\begin{equation}
|z(t,t_0,z^0,\omega(\cdot),\delta (\cdot))|^2 + |\omega(t)|^2
\le R_3^2 -
\frac{t-t_0}{T} (R_3^2-R_1^2)
 \mbox{ for all } \; 
t \in [t_0,t_0+ T]; \label{eq.proof.conditions.for.R_i.8}
\end{equation}
and
\begin{equation*}
\mbox{ if }  \max\limits_{t_0{\le}t{\le}t_0{+}T}|\omega(t)|
\le 3 \sigma_{q{+}2}, \;  \;
|z_0|^2{\le}r_{q{+}2}^2,
  \mbox{ and }  \gamma_k(\parallel \delta(\cdot) {\parallel}_{L_{\infty} [t_0,t_0{+}T]}) \le d_q, \mbox{ then }
\end{equation*}
\begin{equation*}
|z(t,t_0,z^0,\omega(\cdot), \delta(\cdot))|^2 {+} |\omega(t)|^2 \le R_{q{+}2}^2 {-}
\frac{t{-}t_0}{T} (R_{q{+}2}^2{-}R_q^2)
\end{equation*}
\begin{equation}
 \mbox{ for all } \; \; t\in [t_0,t_0{+}T]; \; \; \;
- q_0 - 1 \le q \le 0, \; \; \; {q\in \zz}.
\label{eq.proof.conditions.for.R_i.9}
\end{equation}
and, third,
\begin{equation}
- 2 \varkappa < - \frac{r_{q+2}^2 - R_{q-1}^2}{T} \; \; \; \mbox{
for all } \; -q_0 \le q \le 0, \; \; q \in \zz
\label{eq.additional.cond.for.R.qzero.minus1}
\end{equation}
\begin{equation*}
\frac{\partial V_{k{+}1} (y)}{\partial y} (\psi(t,y,\nu(t,y)) {+} \Delta \phi (t,y))
{=} 2
\langle y, \psi (t,y,\nu(t,y)){+} \Delta \phi (t,y) \rangle < {-}2\varkappa
\end{equation*}
\begin{equation*}
\mbox{ whenever } \; \; t\in [0,T], \; \; R_{{-}q_0{-}1} \le
|y|\le r_{{-}q_0{+}1}, \; y \in \rkplusmpplusl
\end{equation*}
\begin{equation}
\mbox{ and } \gamma_k(\parallel \Delta(\cdot) {\parallel}_{L_{\infty} [t_0,t_0{+}T]})  < d_q
\label{eq.proof.conditions.for.R_i.9_a}
\end{equation}

Then, using Condition C2 of Theorem 3 and the induction over $q
\ge  {-}q_0,$ $q\in \zz,$ if $R_{{-}q_0},$ $R_{{-}q_0{+}1},$ ...,
$R_1,$...,$R_{q{+}1},$ and $\sigma_{{-}q_0},$
$\sigma_{{-}q_0{+}1},$ ..., $\sigma_1,$...,$\sigma_{q{+}1},$ are
already constructed for some $q \ge 2,$ we find $R_{q{+}2} > 0$
and $\sigma_{q{+}2}> 0$ such that
\begin{equation}
r_{q{+}2}< \rho_{q{+}2} < R_{q{+}2} < r_{q{+}3}; \; \; \; \; 0<
\sigma_{q{+}2} \le \sigma_{q{+}1}
 \label{eq.proof.conditions.for.R_i.10}
\end{equation}
and such that for every $z_0 \in \rk,$ every $t_0 \in [0,T],$
every $\omega(\cdot)$ in $C([t_0,t_0{+}T];\rmpplusl)$ and every $\delta(\cdot)$ in $L_{\infty} [t_0,t_0{+}T]$
\begin{equation*}
\mbox{ if } \; \; \max\limits_{t_0{\le}t{\le}t_0{+}T}|\omega(t)|
\le 3 \sigma_{q{+}2},  \; \; |z_0|^2 \le
r_{q{+}2}^2, \; \;
\end{equation*}
\begin{equation*}
 \mbox{ and } \; \;  \gamma_k( \parallel \delta(\cdot) {\parallel}_{L_{\infty} [t_0,t_0{+}T]}) < d_q, \; \; \;  \mbox{ then }
\end{equation*}
\begin{equation*}
|z(t,t_0,z^0,\omega(\cdot),\delta(\cdot))|^2 {+} |\omega(t)|^2 \le R_{q{+}2}^2 {-}
\frac{t{-}t_0}{T} (R_{q{+}2}^2{-}R_q^2)
\end{equation*}
\begin{equation}
 \mbox{ for all } \; \; t \in [t_0,t_0{+}T]; \; \; \; q \ge {-}q_0{-}1,
\; \; {q{\in}\zz}. \label{eq.proof.conditions.for.R_i.11}
\end{equation}

 Define
\begin{equation*}
\Xi_{{-}q_0{+}1}:= {\overline B}_{r_{{-}q_0{+}1}} (0)\; \;
\end{equation*}
\begin{equation}
\mbox{ and } \; \; \Xi_{q{+}1}:={\overline B}_{r_{q{+}1}} (0)
{\setminus} B_{r_q}(0),
  \; \; q \ge -q_0+1, \; q\in \zz;
\label{eq.proof.definition.of.Xi_q.12a}
\end{equation}
\begin{equation*}
P_{{-}q_0{+}1}:= {\Xi}_{{-}q_0{+}1}\cap \left(
\rk{\times}{\overline \Omega}_{\sigma_3}(0) \right)
 \; \mbox{ and }
 \end{equation*}
\begin{equation}
 P_{q{+}1}:= \Xi_{q{+}1} \cap \left(
\rk{\times}{\overline \Omega}_{\sigma_{q{+}4}} (0) \right), \;
q{\ge}{-}q_0{+}1,\;  q{\in}\zz;
\label{eq.proof.definition.of.P_q.12}
\end{equation}
\begin{equation*}
E_{{-}q_0{+}1}:= {\Xi}_{{-}q_0{+}1}\cap \left(
\rk{\times}\left({\overline
\Omega}_{2\sigma_2}(0){\setminus}\Omega_{\sigma_4}(0) \right)
\right) \; \; \mbox{ and }
\end{equation*}
\begin{equation*}
 E_{q{+}1}:= \Xi_{q{+}1} \cap \left(
\rk{\times}\left({\overline \Omega}_{2\sigma_{q{+}2}}
(0){\setminus}\Omega_{\sigma_{q{+}4}}(0) \right)\right),
\end{equation*}
\begin{equation}
 \; \;  q \ge {-}q_0 + 1,\; \;  q \in \zz; \label{eq.proof.definition.of.E_q.13}
\end{equation}
\begin{equation*}
G_{{-}q_0{+}1}:= {\Xi}_{{-}q_0{+}1}\cap \left(
\rk{\times}\left(\rmpplusl{\setminus}{\Omega_{2\sigma_3}(0)}\right)
\right); \; \; \mbox{ and }
\end{equation*}
\begin{equation}
 G_{q{+}1}:= \Xi_{q{+}1} {\setminus} \left(
\rk{\times}{\Omega_{2\sigma_{q{+}2}} (0)} \right), \; \;  q {\ge}
{-}q_0 {+} 1, \; q \in \zz; \label{eq.proof.definition.of.G_q.14}
\end{equation}
\begin{equation*}
K_{q{+}1} {:=} \bigcup\limits_{i{=}{-}q_0{+}1}^{q}\left( E_{i{+}1}
{\cup} P_{i{+}1}\right) \;  \mbox{ and } \; H_{q{+}1} {:=}
\bigcup\limits_{i{=}{-}q_0{+}1}^{q} P_{i{+}1},
\end{equation*}
\begin{equation}
\; \; q \ge {-}q_0+1, \; \;
q{\in}\zz;\label{eq.proof.definition.of.K_q.15}
\end{equation}
Then
\begin{equation*}
H_{q{+}1} \subset K_{q{+}1},\; \; \; q \ge {-}q_0 + 1, \; q \in
\zz.
\end{equation*}
Define
\begin{equation*}
\varepsilon_{q{+}1}:=\min\left\{ \frac{r_{-q_0}}{2} ; \;
\frac{\sigma_{q{+}4}^2}{2}; \;
\min\limits_{-q_0-1{\le}i{\le}q{+}1} \left\{ \frac{r_{i{+}1}  -
R_i}{5}\right\}; \right. \end{equation*}
\begin{equation*}  \min\limits_{-q_0{\le}i{\le}q{+}1}
\left\{\frac{R_{i}  - r_i}{5} \right\};
 \min\limits_{-q_0-1{\le}i{\le}q{+}1} \left\{
\frac{r_{i{+}1}^2  - R_i^2}{5}\right\} ; \end{equation*}
\begin{equation} \left.
\min\limits_{-q_0{\le}i{\le}q{+}1} \left\{\frac{R_{i}^2
 - r_i^2}{5} \right\} \right\}, \; \;  q \ge {-}q_0 + 1, \;
q \in \zz; \label{eq.proof.definition.of.delta_q.16}
\end{equation}
\begin{equation}
m_{q+1}{:=} \max\limits_{\begin{array}{c}  [t,z,z_{k+1}]{\in} [0,T] {\times} K_{q+3} \\ \gamma_k(|\delta|)<d_{q{+}3} \end{array} } (2 |\langle
z, g(t,z,z_{k+1})
{+} \sum\limits_{j{=}1}^{N_k} {\Delta}_j {\varphi}_j (t,z)  \rangle| +1 ) \label{eq.def.of.m.q.plus1}
\end{equation}
and
\begin{equation} D_{q{+}1} {:=} \left\{
\begin{array}{l}
\max\{ \frac{R_{q{+}3}^2}{T}, \; 3m_{q+1} \}  \; \; \;  \; \; \;\; \; \mbox{ if }  q{=}{-}q_0{+}1 \\
\max\{\frac{R_{q{+}3}^2 {-} r_{q{-}2}^2}{T},  \;  3m_{q+1} \}    \mbox{ if }  q {\ge} {-}q_0{+} 2,\\
 \end{array}\right. \;  q \in \zz \label{eq.proof.definition.of.Delta_q.17}
\end{equation}

%Fix any $q \in \zz$ such that $q \ge {-}q_0+1.$ For each
%$M_1{\in}\nn,$ define
%\begin{equation*}
%\Lambda_{q{+}1} (M_1) {:=} \{ [t,z,z_{k{+}1}] {\in} [0,T] {\times}
%(G_{q{+}1} {\cup}  E_{q{+}1})\; | \;  \exists v {\in} \rmpplust
%\end{equation*}
%\begin{equation*}
%(|v| \le  M_1) \wedge  \left( \langle z, g(t,z,z_{k{+}1}) \rangle
%{+} \langle z_{k{+}1}, g_{k{+}1} (t,z,z_{k{+}1},v) \rangle
% {<}\right.
%\end{equation*}
%\begin{equation*}
% \left. {-} 2
%D_{q{+}1}\right) \wedge
% ( \langle z_{k+1}, g_{k+1}
%(t,z,z_{k+1},v) \rangle  < {-} \frac{ 3\sigma^2}{T} ) \}.
%\end{equation*}
%Then $\Lambda_{q{+}1} (M_1) \subset \Lambda_{q{+}1} (M_1{+}1)$ for
%every $M_1 \in \nn,$ and each $\Lambda_{q{+}1} (M_1)$ is an open
%subset of $[0,T] \times (G_{q{+}1} \cup  E_{q{+}1})$ with respect
%to the standard topology of $[0,T] \times (G_{q{+}1} \cup
%E_{q{+}1})$ generated by the norm of $\rr{\times}\rkplusmpplusl.$
%On the other hand, from Assumption A3', it follows that $[0,T]
%\times (G_{q{+}1} \cup  E_{q{+}1})
% \subset \bigcup\limits_{M_1{=}1}^{\infty} \Lambda_{q{+}1}(M_1);$
%Since $[0,T] \times  (G_{q{+}1} \cup  E_{q{+}1})$ is a compact set,
%there exists $M_1(q)\in \nn$ such that $[0,T]\times (G_{q{+}1}\cup
%E_{q{+}1}) \subset \Lambda_{q{+}1} (M_1(q)).$

 Using Assumption A3' and the compactness
of $[0,T] \times  (G_{q{+}1} \cup  E_{q{+}1})$, for
every $q \ge {-}q_0 + 1,$ $q \in \zz,$ one gets the existence of $M_1(q)\in\nn$
such that
\begin{equation*}
\forall [t,z,z_{k{+}1}] {\in} [0,T] {\times} ( G_{q{+}1} {\cup}
E_{q{+}1}) \; 
  \exists v_{t,z,z_{k{+}1}} {\in} \rmpplust  \mbox{ such that } \; |v_{t,z,z_{k{+}1}}|{\le} M_1(q)   \mbox{ and }
\end{equation*}
\begin{equation*}
\mbox{ and } (\forall \Delta {\in} {\bf R}^{N_{k{+}1}} \;  \;  \gamma_k(|\Delta|){<} d_{q{+}1} \Rightarrow \langle z, g(t,z,z_{k{+}1})
{+} \sum\limits_{j{=}1}^{N_k} {\Delta}_j {\varphi}_j (t,z) \rangle
    \end{equation*}
\begin{equation*}
   + \langle z_{k{+}1},
      g_{k{+}1} (t,z,z_{k{+}1},v_{t,z,z_{k{+}1}})
         {+} \sum\limits_{j{=}1}^{N_{k+1}} {\Delta}_j {\varphi}_{k{+}1,j} (t,z,z_{k{+}1}) \rangle {<}  {-} 2
   D_{q{+}1}) \end{equation*}
\begin{equation}
 \mbox{ and } \; \; \; \; \langle z_{k+1}, g_{k+1}
(t,z,z_{k+1},v_{t,z,z_{k{+}1}})
{+} \sum\limits_{j{=}1}^{N_{k+1}} {\Delta}_j {\varphi}_{k{+}1,j} (t,z) \rangle  < {-}\frac{3\sigma^2}{T}
\label{eq.proof.definition.of.M1(q).18}
\end{equation}
In addition,  using Assumption A3' and the compactness
of $[0,T]\times P_{q{+}1}$, for every $q \ge {-}q_0 + 1,$ $q \in
\zz, $ we obtain the existence of $M_2(q)\in\nn$ such that
%%%%%%%%%\begin{equation*}
%%%%%%%%%\forall (t,z,z_{p{+}1}){\in}[0,T]{\times} E_{q{+}1}  \; \exists
%%%%%%%%%v_{t,z,z_{p{+}1}} {\in} \rmpplust{:}  \; \; |v_{t,z,z_{p{+}1}}|{\le}
%%%%%%%%%\end{equation*}
%%%%%%%%%\begin{equation*}
%%%%%%%%%{\le} M_2(q) \; \mbox{ and }  \;  \langle z_{p{+}1}, g_{p{+}1}
%%%%%%%%%(t,z,z_{p{+}1},v_{t,z,z_{p{+}1}}) \rangle {<}
%%%%%%%%%\end{equation*}
%%%%%%%%%\begin{equation}
%%%%%%%%%< - \frac{3\sigma}{T}
%%%%%%%%%   \label{eq.proof.definition.of.M2(q).19}
%%%%%%%%%\end{equation}
\begin{equation*}
{\forall} [t,z,z_{k{+}1}] {\in} [0,T] {\times}  P_{q{+}1} \; \;  
\exists w_{t,z,z_{k{+}1}} {\in} \rmpplust  \mbox{ such that }  \; \;
|w_{t,z,z_{k{+}1}}| \le M_2(q)
\end{equation*}
\begin{equation}
 \; \mbox{ and }
 | \langle z_{k{+}1}, g_{k{+}1}
(t,z,z_{k{+}1},w_{t,z,z_{k{+}1}})
  \rangle | = 0
 \label{eq.proof.definition.of.M3(q).20}
\end{equation}
For each $q\ge{-}q_0+1,$ $q\in\zz, $ define
\begin{equation*}
M(q):= \max \{M_1(q), M_2(q), \max\limits_{{
\begin{array}{l}
 0{\le}t{\le}T\\
 |y|{\le}{2 r_{-q_0}}
 \end{array}}} |\nu(t,y)| \}
\end{equation*}
\begin{equation}
U_q:= \{ u\in\rmpplust \; | \; \;  |u|\le M(q)  \}, \; \; \;
q\ge{-}q_0+1, \;  q\in\zz,
\label{eq.proof.definition.of.M(q).and.U_q.21}
\end{equation}
Without loss of generality, we assume that
\begin{equation*}
M(q) \le M(q{+}1),  \; \; \mbox{ i.e., } \; \;  U_{q} \subset
U_{q{+}1}
\end{equation*}
\begin{equation}
\; \; \mbox{ for all }  \;  q\ge{-}q_0+1, \;  q{\in}\zz.
\label{eq.proof.definition.of.M(q).and.U_q.22}
\end{equation}
Using the compactness of all $U_q,$ ${\overline B}_{r_q}(0),$
take any sequence  ${\{L_q {\}}_{q{=}{-}q_0{+}1}^{\infty}
\subset \rr}$ such that
\begin{equation}
0<L_{q{+}1}\le L_q, \; \; \; \; \; \; \; q \ge -q_0 + 1,\; \;
q{\in}\zz, \label{eq.proof.definition.of.L_q.23}
\end{equation}
\begin{equation*}
2 L_q (|\psi (t,y,u)| {+}|\Delta||\phi(t,y)| {+} 1)  {\le}  1  \; \;  \mbox{ for all } t
{\in}  [0,T],
\end{equation*}
\begin{equation*}
\mbox{ whenever } y \in {\overline B}_{r_{q{+}3}} (0), \; \;   u \in
U_{q{+}3}, \; \;  \gamma_{k+1}(|\Delta|) < d_{q{+}3}, \; \;   \Delta{\in}{\bf R}^{N_{k{+}1}},
\end{equation*}
\begin{equation}
 \; \; \; q \ge {-}q_0+1, \; \; q\in\zz.
\label{eq.proof.definition.of.L_q.24}
\end{equation}
For every $L>0,$ by $\G_L$ denote the system of all the sets given
by
\begin{equation*}
\Gamma_{\Theta(\cdot),\vartheta(\cdot),A_{\Theta},A_{\vartheta}}
:= \{ [s,y] \in  \rr \times  \rkplusmpplusl \; | \; \vartheta
(y)\le s \le
\end{equation*}
\begin{equation*}
 \Theta(y) \} \setminus \left( \{ [s,y] \in \rr \times
\rkplusmpplusl \;  | \;    (s{=} \Theta(y)) \wedge (y {\in}
A_{\Theta}) \} \cup \right. \end{equation*}
\begin{equation*}\left. \{ [s,y] \in \rr\times \rkplusmpplusl\; | \;
(s =  \vartheta(y)) \wedge(y \in A_{\vartheta})\} \right),
\end{equation*}
 where $\Theta (\cdot)$ and $\vartheta(\cdot)$
range over the set of all the functions from class
 $C(\rkplusmpplusl ; [0,T])$ such that
\begin{equation*}
 | \Theta(y_1) {-}\Theta(y_2)| {\le} L |y_1{-}y_2|
\mbox{ and }    |\vartheta(y_1) {-} \vartheta(y_2) | {\le} L
|y_1{-}y_2|
\end{equation*}
\begin{equation}
\mbox{ for all } y_1 \in \rkplusmpplusl,   \; y_2 \in
\rkplusmpplusl, \label{eq.proof.Lipschitz.condition.25}
 \end{equation}
and such that $A_{\Theta} \subset \rkplusmpplusl,$  $A_{\vartheta}
\subset \rkplusmpplusl$ range over the set of all subsets of
$\rkplusmpplusl.$ It is straightforward that for each $L>0,$ $\G_L$ is a semi-ring of sets, i.e.,
 first, $\emptyset \in \G_L;$
second, $\Gamma' \cap \Gamma'' \in \G_L$ for each $\Gamma' \in \G_L,$ and each $\Gamma''\in\G_L;$
  third, for each $\Gamma
\in \G_L,$ and each $\Gamma_1 \in \G_L,$ if $\Gamma_1 \subset
\Gamma,$ then there is a finite sequence
$\{\Gamma_i{\}}_{i{=}2}^l \subset \G_L$ such that $\Gamma =
\bigcup\limits_{j{=}1}^{l}\Gamma_j$ and $\Gamma_i \cap \Gamma_j =
\emptyset,$ whenever $i\not= j,$ $\{i,j\} \subset \{1,2,...,l\}.$

Given $[t,y]{=}[t,z,z_{k{+}1}] {\in} [0,T]{\times}
\left(\rkplusmpplusl {\setminus} B_{r_{{-}q_0{+}1}} (0) \right),$
let $q \ge {-}q_0+1,$ $q\in\zz$ be such that $y\in\Xi_{q{+}1}.$ By
the construction (see
(\ref{eq.proof.definition.of.Xi_q.12a})-(\ref{eq.proof.definition.of.G_q.14})),
we obtain
\begin{equation}
\Xi_{q{+}1}\subset P_{q{+}1} \cup E_{q{+}1}\cup G_{q{+}1}\; \; \;
\mbox{ for all } \; q \ge {-}q_0 + 1,\; q \in \zz
\label{eq.Xi.qplus1.in.EPG.qplis1}
\end{equation}
 Then, the following situations
are possible.

%%%%%%%%%1) $y{\in}E_{q{+}1}.$ Then by
%%%%%%%%%(\ref{eq.proof.definition.of.M2(q).19}), and by the definition of
%%%%%%%%%$U_q$ (see (\ref{eq.proof.definition.of.M(q).and.U_q.21})), there
%%%%%%%%%exists $v_{t,z,z_{p{+}1}} {\in} U_q$ such that inequality
%%%%%%%%%(\ref{eq.proof.definition.of.M2(q).19}) holds. Since $g_{p{+}1}$ is
%%%%%%%%%continuous, there exists a set $T_{t,z,z_{p{+}1}} {\in}
%%%%%%%%%\G_{L_{q{+}2}}$ such that $T_{t,z,z_{p{+}1}} {\subset}
%%%%%%%%%[0,T]{\times}\rkplusmpplusl$ and $T_{t,z,z_{p{+}1}} $ is open in
%%%%%%%%%$[0,T]{\times}\rkplusmpplusl$ w.r.t. the standard topology of
%%%%%%%%%$[0,T]{\times}\rkplusmpplusl$ and such that
%%%%%%%%%\begin{equation}
%%%%%%%%%|y' {-} y''| {<} \delta_{q{+}1} \;   \mbox{ for all } \;
%%%%%%%%%(t',y'){\in} T_{t,z,z_{p{+}1}}, \;   (t'',y''){\in}
%%%%%%%%%T_{t,z,z_{p{+}1}} \label{eq.proof.condition.on.T.26_a}
%%%%%%%%%\end{equation}
%%%%%%%%%and
%%%%%%%%%\begin{equation*}
%%%%%%%%%\langle z_{p{+}1}', g_{p{+}1}
%%%%%%%%%(t',z',z_{p{+}1}',v_{t,z,z_{p{+}1}})\rangle < {-}\frac{3
%%%%%%%%%\sigma^2}{T}
%%%%%%%%%\end{equation*}
%%%%%%%%%\begin{equation}
%%%%%%%%% \mbox{ for all } (t',z',z_{p{+}1}') {\in} T_{t,z,z_{p{+}1}}
%%%%%%%%%\label{eq.proof.condition.on.T.26}
%%%%%%%%%\end{equation}

1) $y \in  (G_{q{+}1} \cup E_{q{+}1}).$ Then, by
(\ref{eq.proof.definition.of.M1(q).18}) and
(\ref{eq.proof.definition.of.M(q).and.U_q.21}), there exist
$v_{t,z,z_{k{+}1}} \in  U_q$ and a set $T_{t,z,z_{k{+}1}} \in
\G_{L_{q{+}2}}$ such that $T_{t,z,z_{k{+}1}} \subset [0,T]\times
\rkplusmpplusl,$ $[t,z,z_{k{+}1}] \in T_{t,z,z_{k{+}1}},$ and
$T_{t,z,z_{k{+}1}} $ is open in $[0,T] \times \rkplusmpplusl$ with
respect to its standard topology and such that
\begin{equation}
|y' {-} y''| {<} \varepsilon_{q{+}1} \;  \mbox{ for all } \; [t',y']
{\in} T_{t,z,z_{k{+}1}}, \;  [t{''},y{''}] {\in} T_{t,z,z_{k{+}1}}
\label{eq.proof.condition.on.Q.27_a}
\end{equation}
and
\begin{equation*}
\langle z', g (t',z',z_{k{+}1}') {+} \sum\limits_{j{=}1}^{N_k} {\Delta}_j {\varphi}_j (t',z') \rangle
+\langle z_{k{+}1}',
g_{k{+}1} (t',z',z_{k{+}1}',v_{t,z,z_{k{+}1}})
\end{equation*}
\begin{equation*}
{+} \sum\limits_{j{=}1}^{N_{k+1}} {\Delta}_j {\varphi}_{k{+}1,j} (t',z',z_{k{+}1}')\rangle <
{-}2 D_{q{+}1} \; \;  \mbox{ and }
\end{equation*}
\begin{equation*}
 \langle z_{k{+}1}',
g_{k{+}1} (t',z',z_{k{+}1}',v_{t,z,z_{k{+}1}})
{+} \sum\limits_{j{=}1}^{N_{k+1}} {\Delta}_j {\varphi}_{k{+}1,j} (t',z',z_{k{+}1}')\rangle < {-}\frac{3
\sigma^2}{T}
\end{equation*}
\begin{equation}
  \mbox{ for all } [t',z',z_{k{+}1}'] \in
T_{t,z,z_{k{+}1}} ,
 \mbox{ whenever } \gamma_{k+1}(|\Delta|) < d_{q{+}1}. \label{eq.proof.condition.on.Q.27}
\end{equation}

2)  $y \in P_{q{+}1}.$ Then, by
(\ref{eq.proof.definition.of.M3(q).20}), 
(\ref{eq.proof.definition.of.M(q).and.U_q.21}), there exist
$w_{t,z,z_{k{+}1}} \in U_q$ and a set $S_{t,z,z_{k{+}1}} \in
\G_{L_{q{+}2}}$ such that $S_{t,z,z_{k{+}1}} \subset [0,T] \times
\rkplusmpplusl,$ $[t,z,z_{k{+}1}] \in S_{t,z,z_{k{+}1}},$ and
$S_{t,z,z_{k{+}1}} $ is open in $[0,T] \times \rkplusmpplusl$ with
respect to its standard topology and such that
\begin{equation}
|y' {-}  y''| {<}  \varepsilon_{q{+}1} \;  \mbox{ for all }  \; [t',y']
{\in} S_{t,z,z_{k{+}1}}, \;  [t'',y''] {\in}  S_{t,z,z_{k{+}1}}
\label{eq.proof.condition.on.S.28_a}
\end{equation}
%and
%\begin{equation*}
%|\langle z_{k{+}1}', g_{k{+}1}
%(t',z',z_{k{+}1}',w_{t,z,z_{k{+}1}})
%\end{equation*}
%\begin{equation*}
%+ \sum\limits_{j{=}1}^{N_{k+1}} {\Delta}_j {\varphi}_{k{+}1,j} (t',z',z_{k{+}1}')\rangle|
%\end{equation*}
%\begin{equation}
% < \frac{1}{4} \min \{
%\varkappa, \frac{\sigma_{q{+}3}^2}{2T} \} \; \mbox{ for all }
%[t',z',z_{k{+}1}'] \in S_{t,z,z_{k{+}1}}
%\label{eq.proof.condition.on.S.28}
%\end{equation}
(If $y \in \Xi_{q{+}1} \cap \Xi_{q{+}2},$ i.e.,
$|y| = r_{q{+}1},$ then we choose $T_{t,z,z_{k{+}1}}$  (or
$S_{t,z,z_{k{+}1}}$) and $v_{t,z,z_{k{+}1}}$ (respectively
$w_{t,z,z_{k{+}1}}$) which correspond to the $\Xi_{q{+}2}.$ Then, by
(\ref{eq.proof.conditions.for.R_i.10}), by
(\ref{eq.proof.definition.of.delta_q.16}), and by
(\ref{eq.proof.definition.of.M(q).and.U_q.22}),
(\ref{eq.proof.definition.of.L_q.23}), inequalities
(\ref{eq.proof.condition.on.Q.27_a})-(\ref{eq.proof.condition.on.S.28_a})
will hold for $q,$ if they hold for $(q+1)$ instead of $q$).

All $[0,T] \times (G_{q{+}1}{\cup}E_{q{+}1})$ and all
$[0,T]\times P_{q{+}1}$ are compact in $[0,T] \times
\rkplusmpplusl,$ therefore by
(\ref{eq.proof.definition.of.Xi_q.12a})-(\ref{eq.proof.definition.of.K_q.15}),
and by
(\ref{eq.proof.condition.on.Q.27_a})-(\ref{eq.proof.condition.on.S.28_a}),
 there exist sequences of sets $\{
T_{t_{\lambda},z_{\lambda},z_{k{+}1}^{\lambda}}
{\}}_{\lambda{=}1}^{\infty} ,$ and $\{ S_{{\hat t}_{\eta},{\hat
z}_{\eta},{\hat z}_{k{+}1}^{\eta}} {\}}_{\eta{=}1}^{\infty} ,$
(along with the corresponding sequences $\{
v_{t_{\lambda},z_{\lambda},z_{k{+}1}^{\lambda}}
{\}}_{\lambda{=}1}^{\infty} ,$ and $\{ w_{{\hat t}_{\eta},{\hat
z}_{\eta},{\hat z}_{k{+}1}^{\eta}} {\}}_{\eta{=}1}^{\infty} $) and
there exist sequences of nonnegative integer numbers $\{
\lambda_{q} {\}}_{q{=}{-}q_0{+}1}^{+\infty},$  and
 $\{ \eta_{q} {\}}_{q{=}{-}q_0{+}1}^{+\infty}$ such that
\begin{equation*}
0=\lambda_{{-}q_0{+}1}<\lambda_q <\lambda_{q{+}1}, \; \; \; \; 0 =
\eta_{{-}q_0{+}1} < \eta_q < \eta_{q{+}1}, \end{equation*}
 \begin{equation} \mbox{ for all } \; \;
q \ge {-}q_0 + 2, \; q \in \zz;
\label{eq.proof.lambda.gamma.eta.29}
\end{equation}
and
\begin{equation*}
[0,T] {\times}  (E_{q{+}1} {\cup} G_{q{+}1})  {\subset}
\bigcup\limits_{\lambda{=}\lambda_q{+}1}^{\lambda_{q{+}1}}
T_{t_{\lambda},z_{\lambda},z_{k{+}1}^{\lambda}}, \; \;
[0,T]{\times} P_{q{+}1} {\subset}\end{equation*}
 \begin{equation}
\bigcup\limits_{\eta{=}\eta_q{+}1}^{\eta_{q{+}1}} S_{{\hat
t}_{\eta},{\hat z}_{\eta},{\hat z}_{k{+}1}^{\eta}}, \; \mbox{ for
all } \;  q {\ge} {-}q_0 {+} 1, \; q {\in} \zz;
\label{eq.proof.definition.of.T.Q.S.30}
\end{equation}
and such that
\begin{equation*}
\left([0,T]{\times} (E_{q{+}1} {\cup} G_{q{+}1})\right) {\cap}
T_{t_{\lambda},z_{\lambda},z_{k{+}1}^{\lambda}} {\not=}
\emptyset,\; \; \; \;  \lambda_q {+}1 {\le} \lambda {\le}
\lambda_{q{+}1}
\end{equation*}
%%%%%%%%%\begin{equation*}
%%%%%%%%% \left([0,T]{\times}
%%%%%%%%%G_{q{+}1}\right) {\cap} Q_{{\bar t}_{\gamma},{\bar z}_{\gamma},{\bar
%%%%%%%%%z}_{p{+}1}^{\gamma}} \not= \emptyset,\; \; \; \gamma_q {+}1
%%%%%%%%%{\le}\gamma {\le} \gamma_{q{+}1};
%%%%%%%%%\end{equation*}
\begin{equation}
 \left([0,T]{\times} P_{q{+}1} \right) \cap S_{{\hat t}_{\eta},{\hat
 z}_{\eta},{\hat
z}_{k{+}1}^{\eta}} \not= \emptyset,\; \; \; \; \; \;\; \; \eta_q
{+}1 \le \eta \le  \eta_{q{+}1};
 \label{eq.proof.definition.of.T.Q.S.30a}
\end{equation}
and
\begin{equation*}
T_{t_{\lambda},z_{\lambda},z_{k{+}1}^{\lambda}} \in \G_{L_{q{+}2}}
\; \; \; \; \; \mbox{ for all } \; \; \lambda_q
 + 1 \le \lambda \le \lambda_{q{+}1}
\end{equation*}
%%%%%%%%%\begin{equation*}
%%%%%%%%%Q_{{\bar t}_{\gamma},{\bar z}_{\gamma},{\bar z}_{p{+}1}^{\gamma}}\in
%%%%%%%%%\G_{L_{q{+}2}} \; \; \; \; \; \mbox{ for all } \; \; \gamma_q
%%%%%%%%%{+}1{\le}\gamma{\le} \gamma_{q{+}1}
%%%%%%%%%\end{equation*}
\begin{equation*}
S_{{\hat t}_{\eta},{\hat z}_{\eta},{\hat z}_{k{+}1}^{\eta}}\in
\G_{L_{q{+}2}} \; \; \; \; \; \mbox{ for all } \; \; \eta_q
{+}1{\le}\eta{\le} \eta_{q{+}1}
\end{equation*}
\begin{equation}
\; \; \; \; \; \; \mbox{ for all } \; \; q \ge {-}q_0 + 1, \; q
\in \zz\label{eq.proof.definition.of.T.Q.S.31}
\end{equation}
Define
\begin{equation*}
T_\lambda :=  T_{t_{\lambda},z_{\lambda},z_{k{+}1}^{\lambda}}, \;
\; \;    S_{\eta} :=  S_{{\hat t}_{\eta},{\hat z}_{\eta},{\hat
z}_{k{+}1}^{\eta}},\; \; v_\lambda :=
v_{t_{\lambda},z_{\lambda},z_{k{+}1}^{\lambda}},
\end{equation*}
 \begin{equation} w_{\eta} :=
w_{{\hat t}_{\eta},{\hat z}_{\eta},{\hat z}_{k{+}1}^{\eta}}\; \;
  \mbox{ for all }     \lambda{\in}\nn, \;
\;    \eta{\in} \nn. \label{eq.proof.definition.of.T.Q.S.32}
\end{equation}

 Since $ \G_{L_q}$ are semirings of sets and $\G_{L_{q{+}1}}
\subset \G_{L_q},$ we use
(\ref{eq.proof.definition.of.delta_q.16}),
(\ref{eq.proof.condition.on.Q.27_a}),
(\ref{eq.proof.condition.on.S.28_a}),
(\ref{eq.proof.definition.of.T.Q.S.31}) and Lemma 2 from
\cite{kolmogorov}, p.40, and the induction over $q \ge
{-}q_0+ 1,$ $q \in \zz,$ and obtain the existence of non-empty
sets $\{\Gamma_l {\}}_{l{=}1}^{\infty} =
\{\Gamma_{\Theta_l(\cdot),\vartheta_l(\cdot),
A_{\Theta_l},A_{\vartheta_l} } {\}}_{l{=}1}^{\infty}$ and a
strictly increasing sequence $\{ l_q {\}}_{q{=}{-}q_0}^{{+}\infty}
\subset \zz$ (with $l_{{-}q_0{+}1} = l_{{-}q_0} = 0$ and with
$l_{q{+}1} > l_q$ for all $q \ge  {-} q_0  + 1$) such that

\begin{enumerate}

\item[(a)] $\Gamma_l \in \G_{L_{q{+}1}},$ and  $|y' -  y''| <
\delta_{q{+}1}$  for all $[t',y'] \in \Gamma_l,$  $[t'',y''] \in
\Gamma_l, $ whenever $l_q{+}1 \le l \le l_{q{+}1},$ $q \ge {-}q_0
+ 1$

\item[(b)] $\Gamma_{l'} \cap \Gamma_{l''} = \emptyset,$ whenever
$l' \not= l''$

\item[(c)] $\bigcup\limits_{l{=}1}^{l_{q{+}1}}\Gamma_l=\left(
\bigcup\limits_{\lambda{=}1}^{\lambda_{q{+}1}} T_{\lambda} \right)
\bigcup \left( \bigcup\limits_{\eta{=}1}^{\eta_{q{+}1}}
S_{\eta}\right)$

\item[(d)] For each $\lambda \in \nn,$ and each $\eta \in \nn,$
there exist (and unique due to (b)) finite sets of natural indices
$C(\lambda),$  $D(\eta),$ such that $T_\lambda =
\bigcup\limits_{l{\in}C(\lambda)} \Gamma_l,$  $S_\eta =
\bigcup\limits_{l{\in}D(\eta)} \Gamma_l,$ and the inequalities
$\lambda_q{+}1 \le \lambda \le \lambda_{q{+}1},$ $\eta_q{+}1 \le
\eta \le \eta_{q{+}1}$ (with $q \ge {-}q_0 + 1,$ $q \in \zz$)
imply respectively: $C(\lambda){\subset}\{l_{q{-}1} + 1,l_{q{-}1}
+ 2,...,l_{q{+}1} \},$ and $D(\eta){\subset}\{l_{q{-}1} +
1,l_{q{-}1} + 2,...,l_{q{+}1} \}$

\end{enumerate}

 From properties (b)-(d) and from
(\ref{eq.proof.definition.of.delta_q.16}),
(\ref{eq.proof.condition.on.Q.27_a}),
(\ref{eq.proof.condition.on.S.28_a}),
 (\ref{eq.proof.definition.of.T.Q.S.32}), we obtain
\begin{equation}
\nn=\left( \bigcup\limits_{\lambda{=}1}^{\infty} C(\lambda)
\right) \bigcup \left( \bigcup\limits_{\eta{=}1}^{\infty}
D(\eta)\right) \label{eq.proof.conditions.on.C.D.E.33}
\end{equation}
 Now we define the
feedback $v(\cdot,\cdot),$ which satisfies conditions (i), (ii) and
(iii) of Theorem 3, as follows

\begin{definition}\label{d1}
Take any $l\in \nn,$ and let $q \ge {-}q_0 + 1,$ $q{\in}\zz$ be
such that $l_q + 1 \le l \le l_{q{+}1}.$ Then, by (c) (and by
(b),(d)) $l \in  \left(
\bigcup\limits_{\lambda{=}1}^{\lambda_{q{+}1}} C(\lambda) \right)
\bigcup\left( \bigcup\limits_{\eta{=}1}^{\eta_{q{+}1}}
D(\eta)\right).$ If $l \in
\bigcup\limits_{\lambda{=}1}^{\lambda_{q{+}1}} C(\lambda)$ then by
the construction (see (b)-(d)) there exists $\lambda(l) \in \nn$
such that $\lambda_q {+}1 \le  \lambda(l) \le  \lambda_{q{+}1}$
and such that $\Gamma_l \subset T_{\lambda(l)},$ and in this case we
define
\begin{equation*}
\chi_l := v_{\lambda(l)}
\end{equation*}
If $l \notin \bigcup\limits_{\lambda{=}1}^{\lambda_{q{+}1}}
C(\lambda),$ then by (\ref{eq.proof.conditions.on.C.D.E.33}), $l
\in \bigcup\limits_{\eta{=}1}^{\eta_{q{+}1}} D(\eta),$
%%%%%%%%%%If
%%%%%%%%%%$l{\in}\bigcup\limits_{p{=}1}^{p_{q{+}1}} D(p),$ then
%%%%%%%%%%by (b)-(d) there exists $p(l){\in}\nn$ such that
%%%%%%%%%%$p_q{+}1{\le}p(l){\le}p_{q{+}1}$ and such that
%%%%%%%%%%$\Gamma_l{\subset}Q_{p(l)}$, and in this case by definition put
%%%%%%%%%%\begin{equation*}
%%%%%%%%%%\chi_l := u_{p(l)};
%%%%%%%%%%\end{equation*}
%%%%%%%%%%otherwise $l{\in}\bigcup\limits_{\eta{=}1}^{\eta_{q{+}1}} E(\eta),$
and, by (b)-(d) there exists $\eta(l) \in \nn$ such that
$\eta_q{+}1 \le \eta(l) \le \eta_{q{+}1}$ and such that $\Gamma_l
\subset S_{\eta(l)},$ and in this case we define
\begin{equation*}
\chi_l := w_{\eta(l)}.
\end{equation*}

\end{definition}

Using the induction over $q \ge {-}q_0 + 1,$ $q \in \zz,$ take a
sequence $\{ h_l {\}}_{l{=}1}^{\infty} \subset ]0,{+}\infty[$ of
positive and small enough numbers. %, more specifically any sequence

For every $l{\in}\nn,$ define $\Gamma_l',$ and $\Gamma_l''$ as
follows
\begin{equation*}
\Gamma_l' :=  \{ [s,y] \in  \rr \times  \rkplusmpplusl \; | \;
\vartheta_l (y) +  \frac{h_l}{2} \le  s  \le  \Theta_l(y) -
\frac{h_l}{2} \} \; \; \; \; \;
\end{equation*}
\begin{equation*}
\Gamma_l''  :=  \{ [s,y]  \in  \rr \times \rkplusmpplusl \; | \;
\vartheta_l (y) +  {h_l} \le  s  \le \Theta_l(y) -  {h_l} \},
\end{equation*}
 \begin{equation}  l \in \nn
\label{eq.proof.definition.of.Gamma'.38}
\end{equation}
(In general, some $\Gamma_l',$ and $\Gamma_l''$ are allowed to be
empty).

 Let $p_l(\cdot,\cdot),$ $l \in\nn,$ be functions of
class $C^{\infty}(\rr \times \rkplusmpplusl; [0,1])$ such that
\begin{equation}
0 \le p_l(t,y)\le 1 \; \; \; \; \; \;  \mbox{ for all } [t,y]
\in \rr \times \rkplusmpplusl
\label{eq.proof.definition.of.gamma_l.39}
\end{equation}
\begin{equation}
p_l(t,y) =  0, \; \; \; \; \; \; \; \; \mbox{ whenever } \;
\; [t,y] \notin \Gamma_l'
\label{eq.proof.definition.of.gamma_l.40}
\end{equation}
\begin{equation}
p_l(t,y) = 1, \; \; \; \; \; \; \mbox{ whenever } \; \; [t,y]
\in \Gamma_l'' \label{eq.proof.definition.of.gamma_l.41}
\end{equation}
Let $p(\cdot)$ be any function of class $C^{\infty}
(\rkplusmpplusl; [0,1])$ such that
\begin{equation}
p(y)=1, \; \; \; \; \; \; \; \mbox{ whenever } \; \; y \in
{\overline B}_{r_{{-}q_0{+}1}}(0)
\label{eq.proof.definition.of.gamma.42}
\end{equation}
and
\begin{equation*}
p(y) = 0, \; \; \; \; \; \;  \mbox{ whenever } \; y \in
\rkplusmpplusl{\setminus} B_{\overline{r}} (0)
\end{equation*}
 \begin{equation} \mbox{ with
some } \; \overline{r} \in ]r_{{-}q_0{+}1}, 2r[
\label{eq.proof.definition.of.gamma.43}
\end{equation}

Define the feedback $v(\cdot,\cdot)$ of class $C^{\mu} ([0,T]
\times \rkplusmpplusl;\rmpplust)$ as follows
\begin{equation*}
v(t,y) :=  \sum\limits_{l{=}1}^{\infty} p_l(t,y) \chi_l  +
p(y) (1 -  \sum\limits_{l{=}1}^{\infty}
p_l(t,y))\nu(t,y)
\end{equation*}
 \begin{equation}
   \mbox{ for all } \;  [t,y] \in [0,T]
\times \rkplusmpplusl \label{eq.proof.definition.of.v.44}
\end{equation}
and extend it smoothly and $T-$ periodically onto  $\rr
\times \rkplusmpplusl$

 Using (\ref{eq.proof.definition.of.gamma_l.40}), (b), and the inclusion $\Gamma_l'' \subset  \Gamma_l' \subset  \Gamma_l,$ $l\in
\nn,$
we obtain that, if $\gamma_l(t,y) \not= 0$ for some $l \in \nn$
and some $[t,y] \in [0,T] \times \rkplusmpplusl,$ then $\gamma_{l'}
(t,y) = 0$ for all $l' \not= l;$ therefore 
$v(\cdot,\cdot)$ given by (\ref{eq.proof.definition.of.v.44}) is
well-defined. 
 Furthermore, by the construction, $v(\cdot,\cdot)$
can be $T-$ periodically smoothly extended onto the whole $\rr
\times \rkplusmpplusl.$ (Indeed, by
(\ref{eq.proof.definition.of.gamma_l.40}), for each $y \in \rkplusmpplusl,$ is $h > 0$ such that $\gamma_l(t,y') = 0,$ $l \in \nn,$ and
therefore $v(t,y')= \gamma(y')\beta(t,y')$ for all $t \in [0,h]
\cup[T - h,T]$ and all $y'$ in some small neighborhood of $y.$ Then the $T$ - periodic extension of $v(\cdot,\cdot)$ given by (\ref{eq.proof.definition.of.v.44}) is of class $C^{\infty}.$

Arguing as in \cite{pavlichkov_ge_2009}, Step 5 it is possible to prove that $\{ h_l {\}}_{l{=}1}^{\infty} \subset ]0,{+}\infty[$ can be chosen so small that this feedback extended $T$ -periodically onto the whole ${\bf R} \times {\bf R}^{k+1}$
is well-defined belongs to $C^{\mu}$ and solves the problem, i.e. globally asymptotically stabilizes system (\ref{eq.backstepping.extended.subsystem.2}).
This completes the proof of Theorem 2 and Theorem 3.

\bibliographystyle{ieeetr}
\bibliography{Ref_main}

\end{document}